\documentclass[12pt,a4paper]{amsart}

\newtheorem{theo+}              {Theorem}           [section]
\newtheorem{prop+}  [theo+]     {Proposition}
\newtheorem{coro+}  [theo+]     {Corollary}
\newtheorem{lemm+}  [theo+]     {Lemma}
\newtheorem{exam+}  [theo+]     {Example}
\newtheorem{rema+}  [theo+]     {Remark}
\newtheorem{defi+}  [theo+]     {Definition}

\newenvironment{theorem}{\begin{theo+}}{\end{theo+}}

\newenvironment{corollary}{\begin{coro+}}{\end{coro+}}
\newenvironment{lemma}{\begin{lemm+}}{\end{lemm+}}
\newenvironment{definition}{\begin{defi+}}{\end{defi+}}

\usepackage{amsthm}
\theoremstyle{plain} \theoremstyle{remark}
\newtheorem{remark}{Remark}
\newtheorem{example}{Example}

\def\E{/\kern-1.0em \equiv }

\title{Some classifications of $\infty$-Harmonic maps between Rienmannian manifolds}
\author{Ze-Ping Wang and Ye-Lin Ou}

\address{School of Math and Computer Science,\newline\indent
Guangxi University for Nationalities,\newline\indent Nanning 530006,
P. R. China.\newline\indent\newline\indent Department of
Mathematics,\newline\indent Texas A $\&$ M
University-Commerce,\newline\indent Commerce, TX 75429
USA.\newline\indent E-mail:yelin$\_$ou@tamu-commerce.edu\;(Ou)}

\begin{document}

\title[Some classifications of $\infty$-Harmonic maps]{Some classifications of
$\infty$-Harmonic maps between Riemannian manifolds}


\maketitle

\section*{Abstract}
\begin{quote}
{\footnotesize $\infty$-Harmonic maps are a generalization of
$\infty$-harmonic functions. They can be viewed as the limiting
cases of p-harmonic maps as p goes to infinity. In this paper, we
give complete classifications of linear and quadratic
$\infty$-harmonic maps from and into a sphere, quadratic
$\infty$-harmonic maps between Euclidean spaces. We describe all
linear and quadratic $\infty$-harmonic maps between Nil and
Euclidean spaces, between Sol and Euclidean spaces. We also study
holomorphic $\infty$-harmonic maps between complex Euclidean
spaces.}
\end{quote}

\section{Introduction}
\indent In this paper, we work in the category of smooth objects so
that all manifolds, vector fields, and maps are assumed
to be smooth unless there is an otherwise statement.\\
\vskip0.3cm \indent The infinity Laplace equation
\begin{equation}
\Delta_{\infty}u:=\frac{1}{2}\langle\nabla\,u, \nabla\left|\nabla
u\right|^{2}\rangle=\sum_{i,j=1}^{m}u_{ij}u_{i}u_{j}=0,
\end{equation}
where $u:\Omega \subset \mathbb{R}^{m} \longrightarrow \mathbb{R}$,
$u_{i}=\frac{\partial u}{\partial x^{i}}$ and
$u_{ij}=\frac{\partial^{2} u}{\partial x^{i}\partial x^{j}}$, was
first discovered and studied by G. Aronsson in his study of
``optimal" Lipschitz extension of functions in the
late 1960s (\cite{Ar1}, \cite{Ar2}). \\

To see why this nonlinear and highly degenerate elliptic PDE has
been so fascinating, we recall that the famous minimal surface
equation can be written as
\begin{equation}\notag
(1+\left|\nabla u\right|^{2})\Delta
u+\sum_{i,j=1}^{m}u_{i}u_{j}u_{ij}=0,
\end{equation}
from which we see that the $\infty$-Laplace equation can be obtained
as {\bf harmonic minimal surface equation} meaning the equation for
harmonic functions with minimal graphs.\\

The solutions of the $\infty$-Laplace equation are called
$\infty$-harmonic functions which have the following
interpretations:
\begin{lemma} (see \cite{Ou1})
Let $u : (M^{m},g) \longrightarrow \mathbb{R}$ be a function. Then
the following conditions are equivalent:
\begin{itemize}
\item[(1)] $u$ is an $\infty$-harmonic function, i.e.,
$\Delta_{\infty}u=0$,
\item[(2)] $u$ is horizontally homothetic;
\item[(3)] $\nabla u$ is perpendicular to $\nabla |\nabla u|^
{2}$;
\item[(4)] ${\rm Hess}_{u}(\nabla u, \nabla u)=0$;
\item[(5)] $|\nabla u|^{2}$ is constant along any integral
curve of $\nabla u$.
\end{itemize}
\end{lemma}
Also, the $\infty$-Laplace equation can be viewed (see \cite{Ar1})
as the formal limit, as $p\rightarrow \infty$, of $p$-Laplace
equation
\begin{equation}\notag
\Delta_{p}\,u:=\left|\nabla\,u\right|^{p-2}\left(\Delta\,u+\frac{p-2}
{\left|\nabla\,u\right|^{2}}\Delta_{\infty}\,u\right)=0.
\end{equation}

Finally, the $\infty$-Laplace equation can be viewed as the
Euler-Lagrange equation of the $L^{\infty}$ variational problem of
minimizing
$$
E_{\infty}(u)={\rm ess\,sup}_{\Omega}\left|{\rm d}\,u\right|
$$
among all Lipschitz continuous functions $u$ with given boundary
values on $\partial \Omega$ (see \cite{ACJ}, \cite{Ba}, and
\cite{BEJ} and the references therein for more detailed
background).\\

Recently, a great deal of research work has been done in the study
of the $\infty$-Laplace equation after the work of Crandall and
Lions (see e.g. \cite{CIL}) on the theory of viscosity solutions for
fully nonlinear problems. Many important results have been achieved
and published in, e.g., \cite{ACJ}, \cite{BB}, \cite{Ba},
\cite{BLW1}, \cite{BLW2}, \cite{BEJ}, \cite{Bh}, \cite{CE},
\cite{CEG}, \cite{CIL}, \cite{CY}, \cite{EG}, \cite{EY}, \cite{J},
\cite{JK}, \cite{JLM1}, \cite{JLM2}, \cite{LM1}, \cite{LM2},
\cite{Ob}. \\

On the other hand, the $\infty$-Laplace equation has been found to
have some very interesting applications in areas such as image
processing (see e.g. \cite{CMS}, \cite{Sa}), mass transfer problems
(see e.g. \cite{EG}), and the study of shape metamorphism (see e.g. \cite{CEPB}).\\

The generalization from harmonic functions to harmonic maps between
Riemannian manifolds was so fruitful that it has not only opened new
fields of study in differential geometry, analysis, and topology but
also brought important applications to many branches in mathematics
and theoretical physics. It would be interesting to study maps
between Riemannian manifolds that generalize $\infty$-harmonic
functions. This was initiated in \cite{OW} where the notion of
$\infty$-harmonic maps between Riemannain manifolds was introduced
as a natural generalization of $\infty$-harmonic functions and as
the limit case of $p$-harmonic maps as $p \rightarrow \infty$.\\

\begin{definition}$($\cite{OW}$)$\label{00d}
A map $\varphi:(M, g)\longrightarrow (N, h)$ between Riemannian
manifolds is called an $\infty$-harmonic map if the gradient of its
energy density belongs to the kernel of its tangent map, i.e.,
$\varphi$ is a solution of the PDEs
\begin{equation}\label{00E}
\tau_{\infty}(\varphi):=\frac{1}{2}{\rm d} \varphi ({\rm grad}
\left|{\rm d} \varphi \right|^{2})=0,
\end{equation}
where $\left|{\rm d} \varphi \right|^{2}={\rm
Trace}_{g}\varphi^{*}h$ is the energy density of $\varphi$.
\end{definition}
A direct computation using local coordinates yields (see also
\cite{OW})
\begin{corollary}\label{C1}
In local coordinates, a map $\varphi:(M, g)\longrightarrow (N, h)$
with\\ $\varphi(x)=(\varphi^{1}(x),\ldots, \varphi^{n}(x))$ is
$\infty$-harmonic if and only if
\begin{equation}\label{00E2}
g({\rm grad}\,\varphi^{\alpha}, {\rm grad}\left|{\rm
d}\varphi\right|^{2})=0,\;\;\;\;\;\;\; \alpha=1, 2,..., n.
\end{equation}
\end{corollary}

Clearly, any $\infty$-harmonic function is an $\infty$-harmonic map
by Definition \ref{00d}. It also follows from the definition that
any map between Riemannian manifolds with constant energy density,
i.e., $\left|{\rm d} \varphi \right|^{2}={\rm
Trace}_{g}\varphi^{*}h={\rm constant}$ is an $\infty$-harmonic map.
Thus, the following important and familiar families are all
$\infty$-harmonic maps:
\begin{itemize}
\item totally geodesic maps,
\item isometric immersions,
\item Riemannian submersions,
\item eigenmaps between spheres.
\end{itemize}
Examples of $\infty$-harmonic maps with nonconstant energy density
include the following classes:
\begin{itemize}
\item projections of multiply warped
products (e.g., the projection of the generalized Kasner
spacetimes),
\item equator maps, and
\item radial projections.
\end{itemize}
We refer the readers to \cite{OW} for details of these and many
other examples and other results including methods of constructing
$\infty$-harmonic maps into Euclidean spaces and spheres,
characterizations of $\infty$-harmonic immersions and submersions,
study of $\infty$-harmonic morphisms which can be characterized as
horizontally homothetic submersions, and the transformation
$\infty$-Laplacians under the the conformal change of metrics.

In this paper, we study the classification of $\infty$-harmonic maps
between certain model spaces. We give complete classifications of
linear and quadratic $\infty$-harmonic maps from and into a sphere,
quadratic $\infty$-harmonic maps between Euclidean spaces. We
describe all linear and quadratic $\infty$-harmonic maps between Nil
and Euclidean spaces and between Sol and Euclidean spaces. We also
study holomorphic $\infty$-harmonic maps complex Euclidean spaces.\\

\section{quadratic $\infty$-harmonic maps between Euclidean spaces}
\vskip0.3cm

As we mentioned in Section 1 that any map with constant energy
density is $\infty$-harmonic. It follows that any affine map
$\varphi: \mathbb{R}^{m}\longrightarrow \mathbb{R}^{n}$ with
$\varphi(X)=AX+b $, where $A$ is an $n\times m$ matrix and
$b\in{R}^{n}$ is a constant, is an $\infty$-harmonic map because of
its constant energy density. Note that there are also globally
defined $\infty$-harmonic maps between Euclidean spaces which are
not affine maps. For example, one can check that $\varphi:
\mathbb{R}^{3}\longrightarrow \mathbb{R}^{2}$ given by
$\varphi(x,y,z)=(\cos x +\cos y +\cos z,\; \sin x +\sin y +\sin z)$
is a map with constant energy density $\left|{\rm d} \varphi
\right|^{2}={\rm Trace}_{g}\varphi^{*}h=3$ and hence an
$\infty$-harmonic maps. In this section, we give a complete
classification of $\infty$-harmonic maps between Euclidean spaces
defined by quadratic polynomials. First, we prove the following
lemma which will be used frequently in this paper.
\begin{lemma}\label{LM1}
Let $A_i,\;i=1, 2, \ldots, n$, be symmetric matrices of $m\times m$.
Then, $(\sum_{j=1}^{n}A_j^2)A_i+A_i(\sum_{j=1}^{n}A_j^2)=0$ for all
$i=1, 2\dots, n$ if and only if $A_i=0$ for $i=1, 2\dots, n$
\end{lemma}
\begin{proof}
Suppose otherwise, i.e., one of $A_i$ is not zero, without loss of
generality, we may assume $A_{1}\neq 0$. Then $rank(A_{1})=K$ with
$1\leq K\leq m $. Without loss of generality, we can choose a
suitable orthogonal matrix $T$ such that ${T^{-1}A_{1}T}$ takes the
diagonal form
\begin{eqnarray}
T^{-1}A_{1}T=\left(\begin{array}{cccc}
\lambda_{1}&0&\cdots&0 \\
0&\lambda_{2}&\cdots&0\\
\vdots&\vdots&\ddots&\vdots\\
0&0&\cdots&\lambda_{m} \\
\end{array}\right)
\end{eqnarray}
where $\lambda_{i_{k}}\neq 0,k=1,2,\ldots,K$ \\Note that
\begin{eqnarray}\label{MAP4}
T^{-1}\sum\limits_{j=1}^{n}(A_{j}^{2})T
=\sum\limits_{j=1}^{n}T^{-1}(A_{j}^{2})T
=\sum\limits_{j=1}^{n}(T^{-1}A_{j}T)^{2}\\\notag
\end{eqnarray}
with each $(T^{-1}A_{j}T)^{2}$ being symmetric matrix. It follows
that
\begin{eqnarray}\label{D03}
&&0=T^{-1}0T=T^{-1}\sum_{j=1}^{n}(A_{j}^{2}A_{1}+A_{1}A_{j}^{2})T\\\notag
&&=\sum_{j=1}^{n}(T^{-1}A_{j}T)^{2}T^{-1}A_{1}T+T^{-1}A_{1}T(T^{-1}A_{j}T)^{2}).\\\notag
\end{eqnarray}
This is impossible because the i-th entry in the main diagonal of
the matrix on the right-hand side of Equation (\ref{D03}) takes the
form
\begin{equation}
2\lambda_i(\lambda_i^2+\sum_{j\geq 2}^{n}|(T^{-1}A_jT)^i|^2),
\end{equation}
where $ (T^{-1}A_jT)^i$ denotes the i-th row vector in
$(T^{-1}A_jT)$, and we know that at least one $\lambda_i$ is not
zero. The contradiction proves the Lemma.
\end{proof}

\begin{theorem}\label{DLA2}
Let $\varphi: \mathbb{R}^{m}\longrightarrow \mathbb{R}^{n}$ be a
quadratic map with $\varphi(X)=(X^{t}A_1X,\ldots,X^{t}A_nX)$, where
$X^{t}=(x^1,\ldots, x^m)\in \mathbb{R}^m$. Then, $\varphi$ is an
$\infty$-harmonic map if and only if $\varphi $
 is a constant map.
\end{theorem}

\begin{proof}
A straightforward computation gives:
\begin{eqnarray}\notag
&&\nabla\varphi^{i}=2X^{t}A_{i},\\\notag &&\left|{\rm
d}\varphi\right|^{2}=\delta^{\alpha
\beta}{\varphi_{\alpha}}^{i}{\varphi_{\beta}}^{j}\delta_{ij}
=\sum_{i=1}^{n}g(\nabla\varphi^{i},\nabla\varphi^{i})\\\notag
&&=\sum_{i=1}^{n}\langle2X^{t}A_{i},2X^{t}A_{i}\rangle
=4\sum_{i=1}^{n}X^{t}A_{i}^{2}X,\;\;{\rm and}\\\notag
&&\nabla\left|{\rm d}\varphi\right|^{2}=
8\sum\limits_{i=1}^{n}X^{t}A_{i}^{2}.\notag
\end{eqnarray}
It follows from Corollary \ref{C1} that $\varphi$ is
$\infty$-harmonic if and only if
\begin{equation}\notag
g(\nabla\,\varphi^{i}, \nabla\left|{\rm
d}\varphi\right|^{2})=0,\;\;\;\;\;\;\; i=1, 2, \ldots, n,
\end{equation}
which is equivalent to
\begin{equation}\label{Q01}
X^{t}A_{i}(\sum_{j=1}^{n}A_{j}^{2})X=0,\; i=1, 2, \ldots, n.
\end{equation}
As the coefficient matrix $A_{i}(\sum_{j=1}^{n}A_{j}^{2})$ of the
quadratic form on the left-hand side of (\ref{Q01}) is not symmetric
in general we can rewrite (\ref{Q01}) as
\begin{equation}
X^{t}(A_{i}(\sum_{j=1}^{n}A_{j}^{2})+(\sum_{j=1}^{n}A_{j}^{2})A_{i})X=0,\;
i=1, 2, \ldots, n.
\end{equation}
 Since $A_{i}(\sum_{j=1}^{n}A_{j}^{2})+(\sum_{j=1}^{n}A_{j}^{2})A_{i}$ is a
symmetric matrix of $m\times m$ we conclude that $\varphi$ is
$\infty$-harmonic if and only if
\begin{equation}\label{N36}
A_{i}(\sum_{j=1}^{n}A_{j}^{2})+(\sum_{j=1}^{n}A_{j}^{2})A_{i}=0,\;\;\;\;\;\;\;
i=1, 2, \ldots, n.
\end{equation}

It follows from this and Lemma \ref{LM1} that $A_{i}=0$ for $i=1, 2,
\ldots, n,$ and hence $\varphi(X)=0$, a
constant map, from which we obtain the Theorem.\\
\end{proof}
\begin{theorem}\label{PE4}
Let $\varphi: \mathbb{R}^{m}\longrightarrow \mathbb{R}^{n}$,
$\varphi(X)=(X^{t}A_{1}X,\ldots,X^{t}A_{n}X)+(AX)^{t}+b$  be a
polynomial map, where $A_i$ is an $m\times m$ symmetric matrix for
$i=1, 2, \ldots, n$, $A$ an ${n\times m}$ matrix, and $b
\in\mathbb{R}^{n}$ .
 Then, $\varphi$ is an $\infty$-harmonic map  if and only if $\varphi$ is
an affine map with $\varphi(X)=(AX)^{t}+b $.
\end{theorem}
\begin{proof}
Let $\alpha_{i}\in\mathbb{R}^{m},\;\;i=1,2,\ldots,n$, denote the
i-th row vector of the matrix $A$. Then,
\begin{equation}\label{00E22334}
\begin{array}{lll}
\nabla\varphi^{i}=2X^{t}A_{i}+\alpha_{i},\;i=1, 2, ..., n,
\end{array}
\end{equation}
\begin{equation}\notag
\begin{array}{lll}
\left|{\rm d}\varphi\right|^{2}=g^{\alpha
\beta}{\varphi_{\alpha}}^{i}{\varphi_{\beta}}^{j}\delta_{ij}\\
=\sum\limits_{i=1}^{n}g(\nabla\varphi^{i},\nabla\varphi^{i})\\
=\sum\limits_{i=1}^{n}\langle 2X^{t}A_{i}+\alpha_{i},2X^{t}A_{i}+\alpha_{i}\rangle\\
=4\sum\limits_{i=1}^{n}X^{t}A_{i}^{2}X+\sum\limits_{i=1}^{n}g\langle\alpha_{i},\alpha_{i}\rangle
+4\sum\limits_{i=1}^{n}\alpha_{i}A_{i}X,
\end{array}
\end{equation}
and
\begin{equation}\label{CD09}
\begin{array}{lll}
\nabla\left|{\rm d}\varphi\right|^{2}=
8\sum\limits_{i=1}^{n}X^{t}A_{i}^{2}+4\sum\limits_{i=1}^{n}\alpha_{i}A_{i}.
\end{array}
\end{equation}
Substituting (\ref{00E22334}) and (\ref{CD09}) the $\infty$-harmonic
map equation (\ref{00E2}) we conclude that $\varphi$ is
$\infty$-harmonic if and only if
\begin{equation}\label{CD11}
\begin{array}{lll}
0=g(\nabla\,\varphi^{i}, \nabla\left|{\rm d}\varphi\right|^{2})\\
=\langle 2X^{t}A_{i}+\alpha_{i},8\sum\limits_{j=1}^{n}X^{t}A_{j}^{2}
+4\sum\limits_{j=1}^{n}\alpha_{j}A_{j}\rangle\\
=16\sum\limits_{j=1}^{n}X^{t}A_{i}A_{j}^{2}X
+8\sum\limits_{j=1}^{n}X^{t}A_{i}A_{j}(\alpha_{j})^{t}
+8\sum\limits_{j=1}^{n}(\alpha_{i}A_{j}^{2})X
+4\sum\limits_{j=1}^{n}\alpha_{i}A_{j}(\alpha_{j})^{t}.\\
\end{array}
\end{equation}
Since Equation (\ref{CD11}) if true for arbitrary $X$, it is
actually an identity of polynomial in $X$. By comparing the
coefficients of the leading terms of the polynomials of at both
sides we have that, if $\varphi$ is $\infty$-harmonic, then
\begin{equation}\label{N311}
16X^{t}A_{i}\sum_{j=1}^{n}A_{j}^{2}X=0,\;\;\; i=1, 2, \ldots, n,
\end{equation}
which is the same as Equation (\ref{Q01}). Now we can use Lemma
\ref{LM1} to conclude that if $\varphi$ is $\infty$-harmonic, then
$A_i=0$ for $i=1, 2, \ldots, n$ and hence $\varphi(X)=(AX)^{t}+b$ is
an affine map. The converse statement clearly true because an affine
map has constant energy density. Therefore, we obtain the theorem.
\end{proof}
\begin{remark}
(A) It would be interesting to know if there is any
$\infty$-harmonic maps $\varphi: \mathbb{R}^{m}\longrightarrow
\mathbb{R}^{n}$ defined
by homogeneous polynomials of degree greater than $2$.\\

(B) We also remark that the situation for the $\infty$-harmonic maps
between semi-Euclidean spaces is quite different in that there are
many examples of non-constant $\infty$-harmonic maps between
semi-Euclidean spaces defined by quadratic polynomials, for example,
let $\mathbb{R}^{2}_{1}$ denote the 2-dimensional semi-Euclidean
space with semi-Euclidean metric $ds^2=-dx^2+dy^2$, then one can
check that the quadratic map $\varphi:
\mathbb{R}^{2}_{1}\longrightarrow \mathbb{R}^{2}_{1}$ defined by
$\varphi(x, y)=(12x^2+12y^2, 13x^2+10xy+13y^2)$ is an map with
energy density $\left|{\rm d} \varphi \right|^{2}={\rm
Trace}_{g}\varphi^{*}h=(\varphi^1_1)^2-(\varphi^1_2)^2-(\varphi^2_1)^2+(\varphi^2_2)^2=0$,
hence it is an $\infty$-harmonic map. For more examples and study of
$\infty$-harmonic maps between Semi-Euclidean spaces see \cite{Zh}.
\end{remark}

\section{Linear $\infty$-harmonic maps from and into a sphere}

In this section, we first derive an equation for linear
$\infty$-harmonic map between conformally flat spaces. We then use
it to give a complete classification of linear $\infty$-harmonic
maps between a Euclidean space and a sphere.
\begin{lemma}\label{PE010}
Let $\varphi:(\mathbb{R}^{m}, g=F^{-2}\delta_{ij})\longrightarrow
(\mathbb{R}^{n}, h=\lambda^{-2}\delta_{\alpha \beta})$  with
\begin{equation*}
\varphi (X)=AX=\left(A^{1}X, \cdots, A^{n}X \right),
\end{equation*}
where $A^i$ is the i-th row vector of A, be a linear map between
conformally flat spaces. Then, $\varphi$ is $\infty$-harmonic if and
only if $A=0$, i.e., $\varphi(X)=AX=0$ is a constant map, or
\begin{equation}
\langle
A^\alpha,\nabla(\frac{F}{\lambda\circ\varphi})\rangle=0,\;\;\;\;\alpha=1,2,\ldots,n.
 \end{equation}
where $\langle,\rangle$ is the Euclidean inner product and $\nabla
f$ denotes the gradient of $f$ taken with respect to the Euclidean
metric on $\mathbb{R}^m$.
\end{lemma}
\begin{proof}
 It is easy to check that for the linear map $\varphi:(\mathbb{R}^{m},
g=F^{-2}\delta_{ij})\longrightarrow (\mathbb{R}^{n},
h=\lambda^{-2}\delta_{\alpha \beta})$  with $ \varphi
(X)=AX=\left(A^{1}X, \cdots, A^{n}X\right)$ we have:

\begin{eqnarray*}
\left|{\rm d}\varphi\right|^{2}=
F^{2}\delta^{ij}{\varphi^{\alpha}}_{i}{\varphi^{\beta}}_{j}
(\lambda^{-2}\delta_{\alpha \beta})\circ \varphi=
(\frac{F}{\lambda\circ
\varphi})^2\sum_{i=1}^{n}\sum_{j=1}^{m}a_{ij}^{2}=(\frac{F}{\lambda\circ
\varphi})^2|A|^2.\\\notag
\end{eqnarray*}
By Corollary \ref{C1}, $\varphi$ is $\infty$-harmonic if and only if
\begin{eqnarray}\notag
&&g({\rm grad}\,\varphi^{\alpha}, {\rm grad}\left|{\rm
d}\varphi\right|^{2})=g^{ij}\varphi^{\alpha}_i(\left|{\rm
d}\varphi\right|^{2})_j\\\notag
&&=F^2\delta^{ij}\varphi^{\alpha}_i(\left|{\rm
d}\varphi\right|^{2})_j=F^2\langle A^\alpha,
\nabla(\frac{F}{\lambda\circ \varphi})^2|A|^2\rangle\\\notag
&&=\frac{2F^3|A|^2}{\lambda\circ \varphi}\langle A^\alpha,
\nabla(\frac{F}{\lambda\circ \varphi})\rangle =0,\;\;\;\;\;\;\;
\alpha=1, 2,..., n,
\end{eqnarray}
from which the Lemma follows.
\end{proof}

 Let $(S^{n}, g_{can})$ be the $n$-dimensional sphere with the standard metric. It
is well known that we can identify $(S^{n}\setminus \{N\}, g_{can})$
with $(\mathbb{R}^{n}, \lambda^{-2}\delta_{ij})$, where
$\lambda=\frac{1+\left|x\right|^{2}}{2}$.\\
Using coordinate $\{x_{i}\}$ we can write the components of $g_{U}$
as:
$${\bar g}_{ij}= \lambda^{-2}{\delta}_{ij},\;\; {\bar
g}^{ij}=\lambda^{2}{\delta}_{ij}.$$

 As an application of Lemma \ref{PE010} we give the following
classification of $\infty$-harmonic maps between spheres.

\begin{theorem}
A linear map $\varphi: (\mathbb{R}^{m},F^{-2}\delta_{ij})\equiv
(S^{m}\setminus \{N\}, g_{can})\longrightarrow (\mathbb{R}^{n},
\lambda^{-2}\delta_{ij})\equiv (S^{n}\setminus \{N\}, g_{can})$
between two spheres with $\varphi (X)=\left(A^{1}X, \cdots, A^{n}X
\right)$ is $\infty$-harmonic if and only if $A=0,$ i.e., $\varphi$
is a constant map, or, $A^tA=I_{m\times m},$ i.e., $\varphi$ is an
isometric immersion.
\end{theorem}

\begin{proof}

To prove the theorem, we applying Lemma \ref{PE010} with
$F=\frac{1+\left|X\right|^{2}}{2}$ and
$\lambda=\frac{1+\left|Y\right|^{2}}{2}$ we conclude that $\varphi$
is $\infty$-harmonic if and only if $A=0$, $\varphi(X)=AX=0$ is a
constant map, or
\begin{eqnarray}\notag
&&\langle
A^\alpha,\nabla(\frac{F}{\lambda\circ\varphi})\rangle=\frac{1}{(\lambda\circ\varphi)^2}\langle
A^\alpha,(\lambda\circ\varphi)\nabla
F-F\nabla(\lambda\circ\varphi)\rangle\\\notag
&&=\frac{1}{(\lambda\circ\varphi)^2}\langle
A^\alpha,(\lambda\circ\varphi)\nabla
[\frac{1}{2}(1+|X|^2)]-F\nabla[\frac{1}{2}(1+|AX|^2)]\rangle\\\notag
&&=\frac{1}{2(\lambda\circ\varphi)^2}\langle A^\alpha,(1+|AX|^2)
X-(1+|X|^2)\nabla(X^tA^tAX)\rangle\\\notag
&&=\frac{1}{2(\lambda\circ\varphi)^2}\langle A^\alpha,(1+|AX|^2)
X-(1+|X|^2)A^tAX\rangle =0,\;\;\;\;\alpha=1,2,\ldots,n,
\end{eqnarray}
which is equivalent to
\begin{equation}\label{Xq1}
(1+|AX|^2)AX- (1+|X|^2)AA^tAX=0,
\end{equation}
for any $X\in \mathbb{R}^m$. It follows that Equation (\ref{Xq1}) is
an identity of polynomials. By comparing coefficients we have
\begin{equation}\label{XL1}
\left\{\begin{array}{rl} AX- AA^tAX=0\\ |AX|^2AX- |X|^2AA^tAX=0.
\end{array}\right.
\end{equation}
for any $X\in \mathbb{R}^m$. It is easy to see that Equation
(\ref{XL1}) implies that $A=0$, or,   $A^tA=I_{m\times m}$ and
$|\varphi(X)|^2=|AX|^2=|X|^2$, from which we obtain the theorem.
\end{proof}
For linear maps between a Euclidean space and a sphere we have
\begin{theorem}
$(1)$ A linear map $\varphi: \mathbb{R}^{m}\longrightarrow
(\mathbb{R}^{n}, \lambda^{-2}\delta_{ij})\equiv (S^{n}\setminus
\{N\}, g_{can})$ from a Euclidean space into a sphere
 with $\varphi (X)=\left(A^{1}X, \cdots, A^{n}X \right)$
is $\infty$-harmonic if and only if $A=0,$ i.e., $\varphi$ is a
constant map.\\
$(2)$ A linear map $\varphi: (\mathbb{R}^{m},
\lambda^{-2}\delta_{ij})\equiv (S^{m}\setminus \{N\},
g_{can})\longrightarrow \mathbb{R}^{n}$ from a sphere into a
Euclidean space with $\varphi (X)=\left(A^{1}X, \cdots, A^{n}X
\right)$ is $\infty$-harmonic if and only if $A=0,$ i.e., $\varphi$
is a constant map.
\end{theorem}

\begin{proof}

To prove the first Statement, we applying Lemma \ref{PE010} with
$F=1$ and $\lambda=\frac{1+\left|Y\right|^{2}}{2}$ we conclude that
$\varphi$ is $\infty$-harmonic if and only if $A=0$,
$\varphi(X)=AX=0$ is a constant map, or
\begin{eqnarray}\notag
&&\langle
A^\alpha,\nabla(\frac{1}{\lambda\circ\varphi})\rangle=-\frac{1}{(\lambda\circ\varphi)^2}\langle
A^\alpha,\nabla(\lambda\circ\varphi)\rangle\\\notag
&&=-\frac{1}{(\lambda\circ\varphi)^2}\langle
A^\alpha,\nabla[\frac{1}{2}(1+|AX|^2)]\rangle\\\notag
&&=-\frac{1}{2(\lambda\circ\varphi)^2}\langle
A^\alpha,\nabla(X^tA^tAX)\rangle\\\notag
&&=-\frac{1}{(\lambda\circ\varphi)^2}\langle A^\alpha,A^tAX\rangle
=0,\;\;\;\;\alpha=1,2,\ldots,n,
\end{eqnarray}
which is equivalent to
\begin{equation}\label{X1}
 AA^tAX=0.
\end{equation}
for any $X\in \mathbb{R}^m$. By letting $X=(A^i)^t,\;i=1, \ldots, n$
in Equation (\ref{X1}) we conclude that $\varphi$ is
$\infty$-harmonic if and only if $AA^tAA^t=0$. Note that
$AA^tAA^t=(AA^t)(AA^t)^t=0$ implies that ${\rm Trace}
(AA^t)=\sum\limits_{i=1}^{n} |A^{i}|^{2}=0$. It follows that
$|A|=0$, i.e., $\varphi$ is a constant map. This gives the first
Statement of the theorem.\\

For the second Statement, we apply Lemma \ref{PE010} with
$\lambda=1$ and $F=\frac{1+\left|X\right|^{2}}{2}$ to conclude that
$\varphi$ is $\infty$-harmonic if and only if $A=0$,
$\varphi(X)=AX=0$ is a constant map, or
\begin{eqnarray}\label{S2}
\langle A^\alpha,\nabla F\rangle=\langle A^\alpha,X\rangle=0
\end{eqnarray}
for $\alpha=1, 2,\ldots,n$ and for all $X\in \mathbb{R}^m$. It is
easy to see that Equation (\ref{S2}) implies that $A^\alpha=0$ for
for $\alpha=1, 2,\ldots,n$ and hence $A=0$, i.e., $\varphi$ is a
constant. This completes the proof of the Theorem.
\end{proof}

\section{Quadratic $\infty$-harmonic maps from and into a sphere}

Again, we identify $(S^{n}\setminus \{N\}, g_{can})$ with
$(\mathbb{R}^{n}, \lambda^{-2}\delta_{ij})$, where
$\lambda=\frac{1+\left|X\right|^{2}}{2}$.\\

\begin{theorem}\label{PE3401}
$(1)$ A quadratic map $\varphi: \mathbb{R}^{m}\longrightarrow
(\mathbb{R}^{n}, \lambda^{-2}\delta_{ij})\equiv (S^{n}\setminus
\{N\}, g_{can})$ into sphere with
$\varphi(X)=(X^{t}A_{1}X,X^{t}A_{2}X,\ldots,X^{t}A_{n}X)$ is
$\infty$-harmonic if and only it is a constant map.\\
$(2)$ A quadratic map from a sphere into a Euclidean space\\
$\varphi:(\mathbb{R}^{m},\lambda^{-2}\delta_{ij})\longrightarrow
\mathbb{R}^{n}$ with
$\varphi(X)=(X^{t}A_{1}X,X^{t}A_{2}X,\ldots,X^{t}A_{n}X)$ is
$\infty$-harmonic if and only  it is a constant map.
\end{theorem}
\begin{proof}
For the Statement (1), we compute:
\begin{equation}\label{P3416}
\begin{array}{lll}
\nabla\varphi^{\alpha} = 2X^{t}A_{\alpha},\\\notag \left|{\rm
d}\varphi\right|^{2}=\delta^{\alpha
\beta}{\varphi_{\alpha}}^{i}{\varphi_{\beta}}^{j}\delta_{ij}(\lambda\circ\varphi)^{-2}
=4\sigma^{2}\sum\limits_{j=1}^{n}X^{t}A_{j}^{2}X,
\end{array}
\end{equation}
where $\sigma=\frac{1}{\lambda\circ\varphi}$. A further computation
gives
\begin{equation}\label{P3439}
\begin{array}{lll}
\nabla\left|{\rm d}\varphi\right|^{2} =8\sigma
(X^{t}\sum\limits_{j=1}^{n}A_{j}^{2}X)\nabla \sigma +8\sigma
^{2}X^{t}\sum\limits_{j=1}^{n}A_{j}^{2}.
\end{array}
\end{equation}
The $\infty$-harmonic map equation for $\varphi$ reads
\begin{equation}\label{qS31}
\begin{array}{lll}
0 =16\sigma(X^{t}\sum\limits_{j=1}^{n}A_{j}^{2}X)\langle
X^{t}A_{\alpha},\nabla
\sigma\rangle+16\sigma^2X^{t}A_{\alpha}\sum\limits_{j=1}^{n}A_{j}^{2}X.
\end{array}
\end{equation}
A direct computation yields  $\nabla
\sigma=-2\sigma^{2}(X^tA_{1}y_{1}+X^tA_{2}y_{2}+\ldots+X^tA_{n}y_{n})$.\\
where $y_{\alpha}=X^{t}A_{\alpha}X$ and for $\alpha=1,2,\ldots,n$.
Substituting this into Equation (\ref{qS31}) we have
\begin{equation}\label{qS1108}
\begin{array}{lll}
0 =&&-32\sigma^{3}(X^{t}\sum\limits_{j=1}^{n}A_{j}^{2}X)(
X^{t}A_{\alpha}A_{1}Xy_{1}+\ldots+X^{t}A_{\alpha}A_{n}Xy_{n})\\
&&+ 16\sigma^{2}X^{t}A_{\alpha}\sum\limits_{j=1}^{n}A_{j}^{2}X
\end{array}
\end{equation}
which is equivalent to
\begin{equation}\label{qS1115}
\begin{array}{lll}
0 =-32(X^{t}\sum\limits_{j=1}^{n}A_{j}^{2}X)(
X^{t}A_{\alpha}A_{1}Xy_{1}+\ldots+X^{t}A_{\alpha}A_{n}Xy_{n})\\+\frac{16}{\sigma}X^{t}A_{\alpha}
\sum\limits_{j=1}^{n}A_{j}^{2}X
\end{array}
\end{equation}
i.e.,
\begin{equation}\label{qS1122}
\begin{array}{lll}
0 =-P(X)+8X^{t}A_{\alpha}\sum\limits_{j=1}^{n}A_{j}^{2}X
\end{array}
\end{equation}
where $P(X)$ denotes a polynomial in $X$ of degree greater than $2$.
Noting that the equation is an identity of polynomials we conclude
that if $\varphi$ is $\infty$-harmonic, then
\begin{equation}\label{PL3388}
X^{t}A_{\alpha}\sum\limits_{j=1}^{3}A_{j}^{2}X=0,\;\;\;\;\;\alpha=1,2,\ldots,n,
\end{equation}
 which is exactly the Equation (\ref{Q01}) and
the same arguments used in the proof of Theorem \ref{DLA2} apply to
give the required results.\\

To prove the second statement, let $\nabla f=(f_{1},\ldots, f_{m})$
denotes the Euclidean gradient of function $f$. Then, a
straightforward computation gives:
\begin{equation}\label{P3416}
\begin{array}{lll}
\nabla\varphi^{\alpha} = 2X^{t}A_{\alpha},\\\notag \left|{\rm
d}\varphi\right|^{2}=\lambda^{2}\delta^{\alpha
\beta}{\varphi_{\alpha}}^{i}{\varphi_{\beta}}^{j}\delta_{ij}
=4\lambda^{2}\sum\limits_{j=1}^{n}X^{t}A_{j}^{2}X,
\end{array}
\end{equation}
and
\begin{equation}\label{P3439}
\begin{array}{lll}
\nabla\left|{\rm d}\varphi\right|^{2} =8\lambda
(X^{t}\sum\limits_{j=1}^{n}A_{j}^{2}X)X
+8\lambda^{2}X^{t}\sum\limits_{j=1}^{n}A_{j}^{2},
\end{array}
\end{equation}
where we have used the fact that $\nabla\lambda=X$. It follows that
$\varphi$ is $\infty$-harmonic if and only if
\begin{equation*}
g({\rm grad}\,\varphi^{i}, {\rm grad}\,\left|{\rm
d}\varphi\right|^{2})=0,\;\;\;\;\;\;\; i=1,2,\ldots,n,
\end{equation*}
which is equivalent to
\begin{equation}\label{qS3}
\begin{array}{lll}
0=g^{ij}\varphi^{\alpha}_{i}\left|{\rm d}\varphi\right|^{2}_{j}=
\lambda^2\langle\nabla\varphi^{\alpha}, \nabla\left|{\rm d}\varphi\right|^{2}\rangle\\
=\lambda^2\langle 2X^{t}A_{\alpha},8\lambda
(X^{t}\sum\limits_{j=1}^{n}A_{j}^{2}X)X
+8\lambda^{2}X^{t}\sum\limits_{j=1}^{n}A_{j}^{2}\rangle\\
=16\lambda^3(X^{t}\sum\limits_{j=1}^{n}A_{j}^{2}X)X^{t}A_{\alpha}X+
16\lambda^4X^{t}A_{\alpha}\sum\limits_{j=1}^{n}A_{j}^{2}X
\end{array}
\end{equation}
for all $X\in \mathbb{R}^m$ and for $\alpha=1,2,\ldots,n$. Note that
Equation (\ref{qS3}) is an identity of polynomials since $\lambda$
is also a polynomial. Comparing the coefficients of the polynomials
we conclude that $\varphi$ is $\infty$-harmonic implies that
$X^{t}A_{\alpha}\sum\limits_{j=1}^{n}A_{j}^{2}X=0$ for
$\alpha=1,2,\ldots,n$, which is exactly the Equation (\ref{Q01}) and
the same arguments used in the proof of Theorem \ref{DLA2} apply to
give the required results.
\end{proof}
\begin{remark}
It is well known that any eigenmap between spheres is of constant
energy density, so any eigenmap is $\infty$-harmonic. It would be
interesting to know if there is any $\infty$-harmonic maps between
spheres which is not an eigenmap.
\end{remark}

\section{Linear $\infty$-harmonic maps from and into Nil space}

In this section we will give a complete classification of linear
$\infty$-harmonic maps between Euclidean spaces and Nil space.

\begin{theorem}\label{PE}
Let $(\mathbb{R}^{3},g_{Nil})$ denote Nil space, where the metric
with respect to the standard coordinates $(x,y,z)$ in
$\mathbb{R}^{3}$ is given by $g_{Nil}={\rm d}x^{2}+{\rm
d}y^{2}+({\rm d}z-x{\rm d}y)^{2}$. Then
\begin{enumerate}
\item[(1)]  A linear function $f:(\mathbb{R}^{3},g_{Nil})\longrightarrow \mathbb{R},\;f(x,y,z)=Ax+By+Cz$ is
an $\infty$-harmonic function if and only if $A=0$ or $C=0$.
\item[(2)] A linear map $\varphi:(\mathbb{R}^{3},g_{Nil})\longrightarrow \mathbb{R}^{n}\; (n\ge 2)$  is
$\infty$-harmonic if and only if $\varphi$ is a composition of the
projection $\pi_{1}:(\mathbb{R}^{3},g_{Nil})\longrightarrow
\mathbb{R}^{2},\;\pi_{1}(x,y,z)=(x,y)$ followed by a linear map
$\mathbb{R}^{2}\longrightarrow \mathbb{R}^{n}$, or, $\varphi$ is a
composition of the projection
$\pi_{2}:(\mathbb{R}^{3},g_{Nil})\longrightarrow
\mathbb{R}^{2},\;\pi_{2}(x,y,z)=(y,z)$ followed by a linear map
$\mathbb{R}^{2}\longrightarrow \mathbb{R}^{n}$.
\end{enumerate}
\end{theorem}
\begin{proof}
For Statement (1), we note that it has been proved in \cite{Ou1}
that a linear function $f:(\mathbb{R}^{3},g_{Nil})\longrightarrow
\mathbb{R},\;f(x,y,z)=Ax+By+Cz$ is an $1$-harmonic if and only if it
is horizontally homothetic which is equivalent to $f$ being
$\infty$-harmonic. It was further shown that this is equivalent to
$A=0$ or $C=0$. To prove Statement (2), one can easily compute the
following components of Nil metric:

\begin{equation}\label{NiM}
\left\{\begin{array}{rl}
g_{11}=1,\; g_{12}=g_{13}=0,\;
g_{22}=1+x^{2},\;g_{23}=-x,\;g_{33}=1;\\
g^{11}=1,\;g^{12}=g^{13}=0,\;g^{22}=1,\;g^{23}=x,\;g^{33}=1+x^{2}.
\end{array}\right.
\end{equation}

Let $\varphi:(\mathbb{R}^{3},g_{Nil})\longrightarrow\mathbb{R}^{n}\;
(n\ge 2)$ be a linear map with
\begin{equation}\label{MAP}
\varphi (X)=\left(\begin{array}{ccc} a_{11} & a_{12} &
a_{13}\\a_{21} & a_{22} & a_{23}\\\ldots & \ldots & \ldots\\a_{n1} &
a_{n2} & a_{n3}
\end{array}\right)\left(\begin{array}{ccc} x\\y\\z\end{array}\right).
\end{equation}
A straightforward computation gives the energy density of $\varphi$
as:
\begin{equation}\label{00E202}
\begin{array}{lll}
\left|{\rm d}\varphi\right|^{2}=g^{\alpha
\beta}{\varphi_{\alpha}}^{i}{\varphi_{\beta}}^{j}\delta_{ij}\\
=\sum\limits_{i=1}^{n}(g^{11}(\frac{\partial \varphi^{i}}{\partial
x_{1}})^{2}+g^{22}(\frac{\partial \varphi^{i}}{\partial
x_{2}})^{2}+g^{33}(\frac{\partial \varphi^{i}}{\partial
x_{3}})^{2}+g^{23}\frac{\partial \varphi^{i}}{\partial
x_{2}}\frac{\partial \varphi^{i}}{\partial
x_{3}}+g^{32}\frac{\partial \varphi^{i}}{\partial
x_{3}}\frac{\partial \varphi^{i}}{\partial
x_{2}})\\
=\sum\limits_{i=1}^{n}a_{i3}^{2}x^{2}+2\sum\limits_{i=1}^{n}a_{i2}a_{i3}x
+\sum\limits_{j=1}^{3}\sum\limits_{i=1}^{n}a_{ij}^{2},\\
\end{array}
\end{equation}
and
\begin{equation}\label{00E203}
\begin{array}{lll}
\frac{\partial \left|{\rm d}\varphi\right|^{2}}{\partial x_{1}}
=\frac{\partial \left|{\rm d}\varphi\right|^{2}}{\partial x}
=2\sum\limits_{i=1}^{n}a_{i3}^{2}x+2\sum\limits_{i=1}^{n}a_{i2}a_{i3},\\
\frac{\partial \left|{\rm d}\varphi\right|^{2}}{\partial
x_{2}}=\frac{\partial \left|{\rm d}\varphi\right|^{2}}{\partial y}=0,\\
\frac{\partial \left|{\rm d}\varphi\right|^{2}}{\partial
x_{3}}=\frac{\partial \left|{\rm d}\varphi\right|^{2}}{\partial z}=0.\\

\end{array}
\end{equation}
It follows from Corollary \ref{C1} and (\ref{00E203}) that $\varphi$
is $\infty$-harmonic if and only if

\begin{equation}\label{N3}
\begin{array}{lll}
g(\nabla\,\varphi^{i}, \nabla\left|{\rm d}\varphi\right|^{2})\\
=g^{\alpha
\beta}{\varphi_{\alpha}}^{i}{\left|{\rm d}\varphi\right|^{2}_{\beta}}\\
=g^{11}\frac{\partial \varphi^{i}}{\partial x_{1}}\frac{\partial
\left|{\rm d}\varphi\right|^{2}}{\partial
x_{1}}+g^{22}\frac{\partial \varphi^{i}}{\partial
x_{2}}\frac{\partial \left|{\rm d}\varphi\right|^{2}}{\partial
x_{2}}+g^{33}\frac{\partial \varphi^{i}}{\partial
x_{3}}\frac{\partial \left|{\rm d}\varphi\right|^{2}}{\partial
x_{3}}+g^{23}\frac{\partial \varphi^{i}}{\partial
x_{2}}\frac{\partial \left|{\rm d}\varphi\right|^{2}}{\partial
x_{3}}+g^{32}\frac{\partial \varphi^{i}}{\partial
x_{3}}\frac{\partial \left|{\rm d}\varphi\right|^{2}}{\partial
x_{2}}\\
=2a_{i1}(\sum\limits_{i=1}^{n}a_{i3}^{2}x+\sum\limits_{i=1}^{n}a_{i2}a_{i3})=0,\;\;\;\;\;\;\;
i=1, 2, \ldots, n.
\end{array}
\end{equation}
Solving Equation (\ref{N3}) we have $a_{i1}=0$, for $i=1, 2, \ldots,
n,$ or $a_{i3}=0$, for $i=1, 2, \ldots, n$, from which we conclude
that the linear map $\varphi:(\mathbb{R}^{3},g_{Nil})\longrightarrow
\mathbb{R}^{n}\; (n\ge 2)$ defined by (\ref{MAP}) is
$\infty$-harmonic if and only if
\begin{equation}\label{MAP1}
\varphi (X)=\left(\begin{array}{ccc} 0 & a_{12} & a_{13}\\0 & a_{22}
& a_{23}\\\ldots& \ldots & \ldots\\0 & a_{n2} & a_{n3}
\end{array}\right)\left(\begin{array}{ccc} x\\y\\z\end{array}\right).
\end{equation}
or
\begin{equation}\label{MAP2}
\varphi (X)=\left(\begin{array}{ccc} a_{11} & a_{12} & 0\\a_{21} &
a_{22} & 0\\\ldots & \ldots & \ldots\\a_{n1} & a_{n2} & 0
\end{array}\right)\left(\begin{array}{ccc} x\\y\\z\end{array}\right).
\end{equation}
Thus, we obtain the theorem.
\end{proof}
\begin{remark}
$(i)$ We remark that in both cases, the maximum possible rank of the
linear $\infty$-harmonic map $\varphi$ is $2$.\\
$(ii)$ Using the energy density formula (\ref{00E202}) we can check
that in case of (\ref{MAP1}) the linear $\infty$-harmonic map
$\varphi$ has non-constant energy density given by a quadratic
polynomial whilst in case of (\ref{MAP2}) the linear
$\infty$-harmonic map $\varphi$ has constant energy density.\\
$(iii)$ It follows from our Theorem that we can choose to have
submersion\\
$\varphi:(\mathbb{R}^{3},g_{Nil})\longrightarrow\mathbb{R}^{2}$ with
\begin{equation}\label{MAP3}
\varphi (X)=\left(\begin{array}{ccc} 0 & a_{12} & a_{13}\\0 & a_{22}
& a_{23}
\end{array}\right)\left(\begin{array}{ccc} x\\y\\z\end{array}\right)
\end{equation}
so that $\varphi$ has non-constant energy density. Clearly,
$\varphi$ cannot be a Riemannian submersion because the energy
density is not constant.
\end{remark}

\begin{theorem}\label{PE01}
A linear map $\varphi:\mathbb{R}^{m}\longrightarrow (\mathbb{R}^{3},
g_{Nil})$ into Nil space is $\infty$-harmonic if and only if
$\varphi$ is a composition of a linear map
$\mathbb{R}^{m}\longrightarrow \mathbb{R}^{2}$ followed by the
inclusion map $\;i_{1}:\mathbb{R}^{2}\longrightarrow
\mathbb{R}^{3},\;i_{1}(y,z)=(0,y,z)$, or, $\varphi$ is a composition
of a linear map $ \mathbb{R}^{m}\longrightarrow \mathbb{R}^{2}$
followed by the inclusion map $\;i_{2}:\mathbb{R}^{2}\longrightarrow
\mathbb{R}^{3},\;i_{2}(x,z)=(x,0,z)$.
\end{theorem}
\begin{proof}
Let $\varphi:\mathbb{R}^{m}\longrightarrow(\mathbb{R}^{3},g_{Nil})$
be a linear map with
\begin{equation}\label{MAP986}
\varphi (X)=\left(\begin{array}{cccc} a_{11} & a_{12} &
\ldots & a_{1m}\\
a_{21} & a_{22} & \ldots&a_{2m}\\
 a_{31} & a_{32} & \ldots &a_{3m}
\end{array}\right)\left(\begin{array}{cccc}x_{1}\\

x_{2}\\

\vdots\\

x_{m}\\

\end{array}\right).
\end{equation}
We can check that the energy density of $\varphi$ is given by:
\begin{equation}\label{NE2}
\begin{array}{lll}
\left|{\rm d}\varphi\right|^{2}=
\delta^{ij}{\varphi_{i}^{\alpha}{\varphi_{j}}^{\beta} g_{\alpha
\beta}}\circ \varphi\\
=\sum\limits_{j=1}^{m}\left(\sum\limits_{\alpha=1}^{3}(\frac{\partial
\varphi^{\alpha}}{\partial x_{j}})^{2}g_{\alpha \alpha}\circ
\varphi+\frac{\partial \varphi^{2}}{\partial x_{j}}\frac{\partial
\varphi^{3}}{\partial x_{j}}g_{23}\circ \varphi+\frac{\partial
\varphi^{3}}{\partial
x_{j}} \frac{\partial \varphi^{2}}{\partial x_{j}}g_{32}\circ \varphi\right)\\
=\sum\limits_{j=1}^{m}\left(a_{1j}^{2}+a_{2j}^{2}(1+x^{2})+a_{3j}^{2}
-a_{2j}a_{3j}x-a_{3j}a_{2j}x\right)\\
=\sum\limits_{j=1}^{m}a_{2j}^{2}x^{2}-2\sum\limits_{j=1}^{m}a_{2j}a_{3j}x
+\sum\limits_{i=1}^{2}\sum\limits_{j=1}^{3}a_{ij}^{2},\;\;\\
\end{array}
\end{equation}
where $x=a_{11}x_{1}+a_{12}x_{2}+\ldots+a_{1m}x_{m}$. Noting that
the domain manifold is Euclidean space we have
\begin{equation}
\nabla\varphi^{i}=(a_{i1},\,a_{i2},\ldots,a_{im}),\;i=1, 2,
3,\\\notag
\end{equation}
and
\begin{eqnarray*}
\nabla\left|{\rm d}\varphi\right|^{2}=\frac{\partial \left|{\rm
d}\varphi\right|^{2}}{\partial x_{i}}\frac{\partial }{\partial
x_{i}}=(2a_{11}(\sum_{j=1}^{m}a_{2j}^{2}x
-\sum_{j=1}^{m}a_{2j}a_{3j}),\ldots,
2a_{1m}(\sum_{j=1}^{m}a_{2j}^{2}x-\sum_{j=1}^{m}a_{2j}a_{3j})). \\
\end{eqnarray*}
By Corollary \ref{C1}, the $\infty$-harmonic map equation for
$\varphi$ becomes
\begin{equation}\label{MAP06}
\langle\nabla\,\varphi^{i}, \nabla\left|{\rm
d}\varphi\right|^{2}\rangle=0,\;\;\;\;\;\;\; i=1, 2, 3,
\end{equation}
which is equivalent to
\begin{equation}\notag
\begin{array}{lll}
2a_{i1}a_{11}(\sum\limits_{j=1}^{m}a_{2j}^{2}x
-\sum\limits_{j=1}^{m}a_{2j}a_{3j})
+2a_{i2}a_{12}(\sum_{j=1}^{m}a_{2j}^{2}x
-\sum\limits_{j=1}^{m}a_{2j}a_{3j})\\+\ldots
+2a_{im}a_{1m}(\sum\limits_{j=1}^{m}a_{2j}^{2}x
-\sum\limits_{j=1}^{m}a_{2j}a_{3j})\\
=2\sum\limits_{k=1}^{m}a_{ik}a_{1k}(\sum\limits_{j=1}^{m}a_{2j}^{2}x
-\sum\limits_{j=1}^{m}a_{2j}a_{3j})=0,\;\;\;\;\;\;\; i=1, 2, 3,
\end{array}
\end{equation}
or
\begin{equation}\label{N032}
\left\{\begin{array}{rl}
2\sum\limits_{k=1}^{m}a_{1k}^{2}(\sum\limits_{j=1}^{m}a_{2j}^{2}x
-\sum\limits_{j=1}^{m}a_{2j}a_{3j})& =0;\\
2\sum\limits_{k=1}^{m}a_{1k}a_{2k}(\sum\limits_{j=1}^{m}a_{2j}^{2}x-
\sum\limits_{j=1}^{m}a_{2j}a_{3j})& =0;\\
2\sum\limits_{k=1}^{m}a_{1k}a_{3k}(\sum\limits_{j=1}^{m}a_{2j}^{2}x-
\sum\limits_{j=1}^{m}a_{2j}a_{3j})& =0.
\end{array}\right.
\end{equation}
Solving (\ref{N032}) we obtain
\begin{equation}\label{N01065}
\left\{\begin{array}{rl}
2\sum\limits_{k=1}^{m}a_{1k}^{2}=0;\\
2\sum\limits_{k=1}^{m}a_{1k}a_{2k}=0;\\
2\sum\limits_{k=1}^{m}a_{1k}a_{3k}=0.\\
\end{array}\right.
\end{equation}
or
\begin{equation}\label{N01074}
 \sum\limits_{j=1}^{m}a_{2j}^{2}x
-\sum\limits_{j=1}^{m}a_{2j}a_{3j} =0.\\
\end{equation}
From Equation (\ref{N01065}) we have
\begin{equation}\label{N01079}
a_{1k}=0,\;{\rm for}\;k=1, 2, \ldots, m.
\end{equation}
Note that if $a_{1k}\ne 0$,  for some $k=1, 2, \ldots, m$, then
$x=a_{11}x_{1}+a_{12}x_{2}+\ldots+a_{1m}x_{m}:
\mathbb{R}^m\longrightarrow \mathbb{R}$ is an onto map and it
follows that Equation (\ref{N01074}) is true for any $x$ as a
polynomial in $x$. Therefore, we have
\begin{equation}\label{N01099}
a_{2k}=0,\;for\;k=1, 2, \ldots, m,
\end{equation}
from which we conclude that the linear map
$\varphi:(\mathbb{R}^{m}\longrightarrow(\mathbb{R}^{3},g_{Nil})$
defined by (\ref{MAP986}) is $\infty$-harmonic if and only if
\begin{equation}\label{MAP08}
\varphi (X)=\left(\begin{array}{cccc}0 & 0 &
\ldots & 0\\
a_{21} & a_{22} & \ldots&a_{2m}\\
a_{31} & a_{32} & \ldots &a_{3m}\\
\end{array}\right)\left(\begin{array}{cccc}x_{1}\\

x_{2}\\

\vdots\\

x_{m}\\

\end{array}\right),
\end{equation}
or
\begin{equation}\label{MAP07}
\varphi (X)=\left(\begin{array}{cccc} a_{11} & a_{12} &
\ldots & a_{1m}\\
0 & 0 & \ldots & 0\\
a_{31} & a_{32} & \ldots &a_{3m}\\
\end{array}\right)\left(\begin{array}{cccc}x_{1}\\

x_{2}\\

\vdots\\

x_{m}\\

\end{array}\right),
\end{equation}

Thus, we obtain the theorem.
\end{proof}
\begin{remark}
$(i)$ We remark that in both cases, the maximum possible rank of the
linear $\infty$-harmonic map $\varphi$ is $2$.\\
$(ii)$ Using the energy density formula (\ref{NE2}) can check that
in both cases of (\ref{MAP08}) and (\ref{MAP07}), the linear
$\infty$-harmonic map $\varphi$ has constant energy density.
\end{remark}

\section{Linear $\infty$-harmonic maps from and into Sol space}

In this section we give a complete classification of linear
$\infty$-harmonic maps between Euclidean spaces and Sol space.
\begin{theorem}\label{PE1}
Let $(\mathbb{R}^{3},g_{Sol})$ denote Sol space, where the metric
with respect to the standard coordinates $(x,y,z)$ in
$\mathbb{R}^{3}$ is given by $g_{Sol}={e^{2z}}{\rm
d}x^{2}+{e^{-2z}}{\rm d}y^{2}+{\rm d}{z}^{2}$. Then
\begin{enumerate}
\item[(1)]  A linear function $f:(\mathbb{R}^{3},g_{Sol})\longrightarrow \mathbb{R},\;f(x,y,z)=Ax+By+Cz$ is
an $\infty$-harmonic function if and only if $C=0$ or $A=B=0$.
\item[(2)] A linear map $\varphi:(\mathbb{R}^{3},g_{Sol})\longrightarrow \mathbb{R}^{n}\; (n\ge 2)$  is
$\infty$-harmonic if and only if $\varphi$ is a composition of the
projection $\pi_{1}:(\mathbb{R}^{3},g_{Sol})\longrightarrow
\mathbb{R}^{2},\;\pi_{1}(x,y,z)=(x,y)$ followed by a linear map
$\mathbb{R}^{2}\longrightarrow \mathbb{R}^{n}$, or, $\varphi$ is a
composition of the projection
$\pi_{2}:(\mathbb{R}^{3},g_{Sol})\longrightarrow
\mathbb{R}^{2},\;\pi_{2}(x,y,z)=(z)$ followed by a linear map
$\mathbb{R}\longrightarrow \mathbb{R}^{n}$.
\end{enumerate}
\end{theorem}
\begin{proof}
The Statement (1) is proved in \cite{Ou1}. To prove Statement (2),
one can easily compute the following components of Sol metric:
\begin{align}\label{Solm}
&g_{11}=e^{2z},\;g_{22}=e^{-2z},\;g_{33}=1,\;{\rm all \; other},\;g_{ij}=0;\\\
&g^{11}=e^{-2z},\;g^{22}=e^{2z},\;g^{33}=1,\;{\rm all \; other}
,\;g^{ij}=0.\notag
\end{align}
Let $\varphi:(\mathbb{R}^{3},g_{Sol})\longrightarrow\mathbb{R}^{n}\;
(n\ge 2)$ be a linear map with
\begin{equation}\label{MAP01}
\varphi (X)=\left(\begin{array}{ccc} a_{11} & a_{12} &
a_{13}\\a_{21} & a_{22} & a_{23}\\\ldots & \ldots & \ldots\\a_{n1} &
a_{n2} & a_{n3}
\end{array}\right)\left(\begin{array}{ccc} x\\y\\z\end{array}\right).
\end{equation}
A straightforward computation gives the energy density of $\varphi$
as:

\begin{equation}\label{00E205}
\begin{array}{lll}
\left|{\rm d}\varphi\right|^{2}= g^{\alpha
\beta}{\varphi_{\alpha}}^{i}{\varphi_{\beta}}^{j}\delta_{ij}\\\notag
=\sum\limits_{i=1}^{n}g^{11}(\frac{\partial \varphi^{i}}{\partial
x_{1}})^{2}+g^{22}(\frac{\partial \varphi^{i}}{\partial
x_{2}})^{2}+g^{33}(\frac{\partial \varphi^{i}}{\partial
x_{3}})^{2}\\
=\sum\limits_{i=1}^{n}a_{i1}^{2}e^{-2z}
+\sum\limits_{i=1}^{n}a_{i2}^{2}e^{2z}+\sum\limits_{i=1}^{n}a_{i3}^{2},\\\notag
\end{array}
\end{equation}
and
\begin{equation}\label{00E206}
\begin{array}{lll}
\frac{\partial \left|{\rm d}\varphi\right|^{2}}{\partial x_{1}}
=\frac{\partial \left|{\rm d}\varphi\right|^{2}}{\partial x}=0,\\
\frac{\partial \left|{\rm d}\varphi\right|^{2}}{\partial x_{2}}
=\frac{\partial \left|{\rm d}\varphi\right|^{2}}{\partial y}=0,\\
\frac{\partial \left|{\rm d}\varphi\right|^{2}}{\partial x_{3}}
=\frac{\partial \left|{\rm d}\varphi\right|^{2}}{\partial z}=
-2\sum\limits_{i=1}^{n}a_{i1}^{2}e^{-2z}+2\sum\limits_{i=1}^{n}a_{i2}^{2}e^{2z}.\\
\end{array}
\end{equation}
It follows from Corollary \ref{C1} and (\ref{00E206}) that $\varphi$
is $\infty$-harmonic if and only if

\begin{equation}\label{N31}
\begin{array}{lll}
g(\nabla\,\varphi^{i}, \nabla\left|{\rm d}\varphi\right|^{2})\\
=g^{\alpha
\beta}{\varphi_{\alpha}}^{i}{\left|{\rm d}\varphi\right|^{2}_{\beta}}\\
=g^{11}\frac{\partial \varphi^{i}}{\partial x_{1}}\frac{\partial
\left|{\rm d}\varphi\right|^{2}}{\partial
x_{1}}+g^{22}\frac{\partial \varphi^{i}}{\partial
x_{2}}\frac{\partial \left|{\rm d}\varphi\right|^{2}}{\partial
x_{2}}+g^{33}\frac{\partial \varphi^{i}}{\partial
x_{3}}\frac{\partial \left|{\rm d}\varphi\right|^{2}}{\partial
x_{3}}\\
=-2a_{i3}(\sum\limits_{i=1}^{n}a_{i1}^{2}e^{-2z}-\sum\limits_{i=1}^{n}a_{i2}e^{2z})=0,\;\;\;\;\;\;\;
i=1, 2, \ldots, n.
\end{array}
\end{equation}
Solving Equation (\ref{N31}) we have $a_{i3}=0$, for $i=1, 2,
\ldots, n,$ or $a_{i1}=a_{i2}=0$, for $ i=1, 2, \ldots, n$, from
which we conclude that the linear map
$\varphi:(\mathbb{R}^{3},g_{Sol})\longrightarrow \mathbb{R}^{n}\;
(n\ge 2)$ defined by (\ref{MAP01}) is $\infty$-harmonic if and only
if
\begin{equation}\label{MAP4}
\varphi (X)=\left(\begin{array}{ccc} 0 & 0 & a_{13}\\0 & 0 &
a_{23}\\\ldots& \ldots & \ldots\\0 & 0
 & a_{n3}
\end{array}\right)\left(\begin{array}{ccc}
x\\y\\z\end{array}\right),
\end{equation}
or
\begin{equation}\label{MAP5}
\varphi (X)=\left(\begin{array}{ccc} a_{11} & a_{12} & 0\\a_{21} &
a_{22} & 0\\\ldots & \ldots & \ldots\\a_{n1} & a_{n2} & 0
\end{array}\right)\left(\begin{array}{ccc} x\\y\\z\end{array}\right).
\end{equation}
Thus, we obtain the theorem.
\end{proof}
\begin{remark}
It follows from Theorem \ref{PE1} that the maximum rank of the
linear $\infty$-harmonic maps from Sol space into Euclidean space is
$2$. In case of (\ref{MAP5}) the linear $\infty$-harmonic map has
non-constant energy density.
\end{remark}
\begin{theorem}\label{tosol}
A linear map $\varphi:\mathbb{R}^{m}\longrightarrow (\mathbb{R}^{3},
g_{Sol})$ is $\infty$-harmonic if and only if it is a composition of
a linear map $\mathbb{R}^{m}\longrightarrow \mathbb{R}^{2}$ followed
by the inclusion map $\;i_{1}:\mathbb{R}^{2}\longrightarrow
\mathbb{R}^{3},\;i_{1}(x,y)=(x,y,0)$, or, $\varphi$ is a composition
of a linear map $\mathbb{R}^{m}\longrightarrow \mathbb{R}$
 followed by an inclusion map $\;i_{2}:\mathbb{R}\longrightarrow
\mathbb{R}^{3},\;i_{2}(z)=(0,0,z)$.
\end{theorem}
\begin{proof}
For the linear map
$\varphi:\mathbb{R}^{m}\longrightarrow(\mathbb{R}^{3},g_{Sol})$
 with
\begin{equation}\label{MAP1164}
\varphi (X)=\left(\begin{array}{cccc} a_{11} & a_{12} &
\ldots & a_{1m}\\
a_{21} & a_{22} & \ldots&a_{2m}\\
a_{31} & a_{32} & \ldots &a_{3m}
\end{array}\right)\left(\begin{array}{cccc}  x_{1}\\

x_{2}\\

\vdots\\

x_{m}\\
\end{array}\right),
\end{equation}
we have:
\begin{equation}
\nabla\varphi^{i}=(a_{i1},\,a_{i2},\ldots,a_{im}),\;i=1, 2,
3,\\\notag
\end{equation}
and the energy density
\begin{equation}
\begin{array}{lll}
\left|{\rm d}\varphi\right|^{2}=
\delta^{ij}{\varphi_{i}^{\alpha}{\varphi_{j}}^{\beta} g_{\alpha
\beta}}\circ \varphi\\
=\sum\limits_{j=1}^{m}\sum\limits_{\alpha=1}^{3}(\frac{\partial
\varphi^{\alpha}}{\partial x_{j}})^{2}g_{\alpha
\alpha}\circ \varphi\\
=\sum\limits_{j=1}^{m}a_{1j}^{2}e^{2z}+\sum\limits_{j=1}^{m}a_{2j}^{2}e^{-2z}
+\sum\limits_{j=1}^{m}a_{3j}^{2},\;\;\\\notag
\end{array}
\end{equation}
where $z=a_{31}x_{1}+a_{32}x_{2}+\ldots+a_{3m}x_{m}$. Since the
domain manifold in a Euclidean space, tt follows from Corollary
\ref{C1} that the $\infty$-harmonic map equation for $\varphi$
becomes
\begin{equation}\label{MAP06}
\langle\nabla\,\varphi^{i}, \nabla\left|{\rm
d}\varphi\right|^{2}\rangle=0,\;\;\;\;\;\;\; i=1, 2, 3,
\end{equation}
which is equivalent to which is equivalent to
\begin{equation}
\begin{array}{lll}
2a_{i1}a_{31}(\sum\limits_{j=1}^{m}a_{1j}^{2}e^{2z}-\sum\limits_{j=1}^{m}a_{2j}^{2}e^{-2z})
+2a_{i2}a_{32}(\sum\limits_{j=1}^{m}a_{1j}^{2}e^{2z}-\sum\limits_{j=1}^{m}a_{2j}^{2}e^{-2z})\\+\ldots
+2a_{im}a_{3m}(\sum\limits_{j=1}^{m}a_{1j}^{2}e^{2z}-\sum\limits_{j=1}^{m}a_{2j}^{2}e^{-2z})\\
=2\sum\limits_{k=1}^{m}a_{ik}a_{3k}(\sum\limits_{j=1}^{m}a_{1j}^{2}e^{2z}
-\sum\limits_{j=1}^{m}a_{2j}^{2}e^{-2z})=0,\;\;\;\;\;\;\; i=1, 2, 3,
\end{array}
\end{equation}
or
\begin{equation}\label{N0032}
\left\{\begin{array}{rl}
2\sum\limits_{k=1}^{m}a_{1k}a_{3k}(\sum\limits_{j=1}^{m}a_{1j}^{2}e^{2z}-\sum\limits_{j=1}^{m}a_{2j}^{2}e^{-2z})& =0\\
2\sum\limits_{k=1}^{m}a_{2k}a_{3k}(\sum\limits_{j=1}^{m}a_{1j}^{2}e^{2z}-\sum\limits_{j=1}^{m}a_{2j}^{2}e^{-2z})& =0\\
2\sum\limits_{k=1}^{m}a_{3k}^{2}(\sum\limits_{j=1}^{m}a_{1j}^{2}e^{2z}-\sum\limits_{j=1}^{m}a_{2j}^{2}e^{-2z})& =0.\\
\end{array}\right.
\end{equation}
From this we have either
\begin{equation}\label{N01232}
\left\{\begin{array}{rl}
2\sum\limits_{k=1}^{m}a_{1k}a_{3k}=0;\\
2\sum\limits_{k=1}^{m}a_{2k}a_{3k}=0;\\
2\sum\limits_{k=1}^{m}a_{3k}^{2}=0,
\end{array}\right.
\end{equation}
or
\begin{equation}\label{N01240}
(\sum\limits_{j=1}^{m}a_{1j}^{2}e^{2z}-\sum\limits_{j=1}^{m}a_{2j}^{2}e^{-2z})
=0.
\end{equation}

Solving Equation (\ref{N01232}) we have
\begin{equation}\label{N01245}
a_{3k}=0,\;for\;k=1, 2, \ldots,m.
\end{equation}

Solving Equation (\ref{N01240}) we have
\begin{equation}\label{N01249}
a_{1k}=a_{2k}=0,\;for\;k=1, 2,\ldots,m.
\end{equation}
The Theorem follows from (\ref{N01245}) and (\ref{N01249}).
\end{proof}
\begin{remark}
It follows from Theorem \ref{tosol} that the maximum rank of the
linear $\infty$-harmonic maps from  Euclidean space into Sol space
is $2$, and any linear $\infty$-harmonic map into Sol space has
constant energy density.
\end{remark}

\section{Quadratic $\infty$-harmonic maps into Sol and Nil spaces}

\begin{theorem}\label{DLA3127}
Let $(\mathbb{R}^{3}, g_{Sol})$ denote Sol space, where the metric
with respect to the standard coordinates $(x,y,z)$ in
$\mathbb{R}^{3}$ is given by $g_{Sol}={e^{2z}}{\rm
d}x^{2}+{e^{-2z}}{\rm d}y^{2}+{\rm d}{z}^{2}$. Then,  a quadratic
map $\varphi: \mathbb{R}^{m}\longrightarrow
(\mathbb{R}^{3},g_{Sol})$ with
$\varphi(X)=(X^{t}A_{1}X,X^{t}A_{2}X,X^{t}A_{3}X)$ is an
$\infty$-harmonic map  if  and  only if it
 is a constant map.
\end{theorem}

\begin{proof}
One can easily compute the following components of Sol metric:
\begin{align}\nonumber
&h_{11}=e^{2z},\;h_{22}=e^{-2z},\;h_{33}=1,\;{\rm all \; other}\;h_{ij}=0;\\\
&h^{11}=e^{-2z},\;h^{22}=e^{2z},\;h^{33}=1,\;{\rm all \; other}
\;h^{ij}=0.\notag
\end{align}
A straightforward computation gives:
\begin{equation}\label{00E23146}
\begin{array}{lll}
\nabla\varphi^{i}=2X^{t}A_{i},\\\notag \left|{\rm
d}\varphi\right|^{2}=\delta^{\alpha
\beta}{\varphi_{\alpha}}^{i}{\varphi_{\beta}}^{j}h_{ij}\circ
\varphi=\langle \nabla \varphi^{i},\nabla \varphi^{j}\rangle
h_{ij}\circ \varphi\\\notag =|\nabla\,\varphi^{i}|^{2}\,h_{ii}\circ
\varphi
=(4X^{t}A_{1}^{2}X)e^{2z}+(4X^{t}A_{2}^{2}X)e^{-2z}+4X^{t}A_{3}^{2}X,
\end{array}
\end{equation}
where $z=X^{t}A_{3}X$.\\

Furthermore,
\begin{equation}\label{00E23260}
\begin{array}{lll}
\nabla\left|{\rm d}\varphi\right|^{2}
= 8X^{t}A_{1}^{2}e^{2z}+8X^{t}A_{1}^{2}X e^{2z}\nabla z\\
+8X^{t}A_{2}^{2}e^{-2z}-8X^{t}A_{2}^{2}X e^{-2z}\nabla z\\
+8X^{t}A_{3}^{2}\\
=8X^{t}A_{1}^{2}e^{2z}+(16X^{t}A_{1}^{2}X e^{2z})X^{t}A_{3}^{2}\\
+8X^{t}A_{2}^{2}e^{-2z}-(16X^{t}A_{2}^{2}X e^{-2z})X^{t}A_{3}^{2}\\
+8X^{t}A_{3}^{2},\\\notag
\end{array}
\end{equation}
and
\begin{equation*}
\begin{array}{lll}
g(\nabla\,\varphi^{i}, \nabla\left|{\rm d}\varphi\right|^{2})\\
=\langle 2X^{t}A_{i}, 8X^{t}A_{1}^{2}e^{2z} +(16X^{t}A_{1}^{2}X
e^{2z})X^{t}A_{3}^{2} +8X^{t}A_{2}^{2}e^{-2z}-(16X^{t}A_{2}^{2}X
e^{-2z})X^{t}A_{3}^{2}
+8X^{t}A_{3}^{2}\rangle\\
=16(X^{t}A_{i}A_{1}^{2}X+2X^{t}A_{1}^{2}X
X^{t}A_{i}A_{3}X)e^{2z}\\
+16(X^{t}A_{i}A_{2}^{2}X-2X^{t}A_{2}^{2}X
X^{t}A_{i}A_{3}X)e^{-2z}\\
+16X^{t}A_{i}A_{3}^{2}X.\\
\end{array}
\end{equation*}
It follows from the $\infty$-harmonic map equation (\ref{00E2}) that
$\varphi$ is $\infty$-harmonic if and only if
\begin{equation}\label{P3207}
\begin{array}{lll}
16(X^{t}A_{i}A_{1}^{2}X+X^{t}A_{1}^{2}X
X^{t}A_{i}A_{3}X)e^{2z}\\
+16(X^{t}A_{i}A_{2}^{2}X-X^{t}A_{2}^{2}X
X^{t}A_{i}A_{3}X)e^{-2z}\\
+16X^{t}A_{i}A_{3}^{2}X=0,\; \;\;\;{\rm for\; all}\,
X\in\mathbb{R}^m,\;{\rm and}\; i=1, 2, 3.
\end{array}
\end{equation}
Note that Equation (\ref{P3207}) is an identity of functions which
are analytic. We can substitute the Taylor expansions for
$e^{2X^tA_{3}X}$ and $e^{-2X^tA_{3}X}$ into (\ref{P3207}) and
compare the coefficients of the second degree terms to get
\begin{equation}\notag
\begin{array}{lll}
16(X^{t}A_{i}\sum\limits_{j=1}^{3}A_{j}^{2}X)
=0,\;\;\;\;\;\;\;\;i=1,2,3.
\end{array}
\end{equation}
From this we obtain
\begin{equation}\notag
A_{i}(\sum\limits_{j=1}^{3}A_{j}^{2})+(\sum\limits_{j=1}^{3}A_{j}^{2})A_{i}=0,
\end{equation}
and Lemma \ref{LM1} applies to complete the proof of the Theorem.
\end{proof}

\begin{theorem}\label{PE3268}
A quadratic map $\varphi: \mathbb{R}^{m}\longrightarrow
(\mathbb{R}^{3},g_{Nil})$,
$\varphi(X)=(X^{t}A_{1}X,X^{t}A_{2}X,X^{t}A_{3}X)$ into Nil space is
$\infty$-harmonic if  and  only if it is a constant map.
\end{theorem}
\begin{proof}
Using the components of Nil metric (\ref{NiM}) with notations
$g_{ij}$ replaced by $h_{ij}$ one can easily check that:
\begin{equation*}
\begin{array}{lll}
\nabla\varphi^{i}=2X^{t}A_{i},\\ \left|{\rm
d}\varphi\right|^{2}=\delta^{\alpha
\beta}{\varphi_{\alpha}}^{i}{\varphi_{\beta}}^{j}h_{ij}\circ\varphi\\
=\sum_{\alpha=1}^{m}{\varphi_{\alpha}}^{i}{\varphi_{\alpha}}^{j}h_{ij}\\
=(4X^{t}A_{2}^{2}X)x^{2}+4\sum\limits_{j=1}^{3}X^{t}A_{j}^{2}X
-4(X^{t}A_{2}A_{3}+A_3A_2X)x,
\end{array}
\end{equation*}
where $x=X^{t}A_{1}X$, and \\
\begin{equation*}
\begin{array}{lll}
\nabla\left|{\rm d}\varphi\right|^{2}=\frac{\partial \left|{\rm
d}\varphi\right|^{2}}{\partial x_{i}}\frac{\partial}{\partial
x_{i}}\\
=(8X^{t}A_{2}^{2})x^{2}+(16X^{t}A_{2}^{2}X)(X^{t}A_{1})x
+8\sum\limits_{j=1}^{3}X^{t}A_{j}^{2}\\
-8(X^{t}A_{2}A_{3}+X^{t}A_{3}A_{2})x-8(X^{t}A_{2}A_{3}+A_3A_2X)X^{t}A_{1}.
\end{array}
\end{equation*}

By Corollary \ref{C1}, the $\infty$-harmonic map equation for
$\varphi$ reads
\begin{equation*}
\begin{array}{lll}
0=\langle\nabla\,\varphi^{i}, \nabla\left|{\rm d}\varphi\right|^{2}\rangle\\
=\langle
2X^{t}A_{i},(8X^{t}A_{2}^{2})x^{2}+(16X^{t}A_{2}^{2}X)(X^{t}A_{1})x
+8\sum\limits_{j=1}^{3}X^{t}A_{j}^{2}
-8(X^{t}A_{2}A_{3}+X^{t}A_{3}A_{2})x\\
-8(X^{t}(A_{2}A_{3}+A_3A_2X)X^{t}A_{1}\rangle\\
=16X^{t}A_{i}\sum\limits_{j=1}^{3}A_{j}^{2}X
+P(X),\;\;\;\;\;i=1,2,3.
\end{array}
\end{equation*}
where $P(X)$ denotes a polynomial in $X$ of degree greater than $2$.
Noting that the equation is an identity of polynomials we conclude
that if $\varphi$ is $\infty$-harmonic, then

\begin{equation}\label{P3388}
X^{t}A_{i}\sum\limits_{j=1}^{3}A_{j}^{2}X=0,\;\;\;\;\;i=1,2,3,
\end{equation}
which is the same as Equation (\ref{Q01}) for $n=3$ so the rest of
the proof is exactly the same as the part in the proof of Theorem
\ref{DLA2}.
\end{proof}

\begin{example}
We remark that there are many polynomial $\infty$-harmonic maps
$\varphi:(\mathbb{R}^{3},g_{Nil})\longrightarrow \mathbb{R}^{m}$,
for instance,
\begin{equation}\notag
\varphi:(\mathbb{R}^{3},g_{Nil})\longrightarrow \mathbb{R}^{2}
\end{equation}
with $\varphi(x,y,z)=(z-xy/2,2z-xy)$ is an $\infty$-harmonic map
which has nonconstant energy density $\left|{\rm d} \varphi
\right|^{2}=5(1+(x^{2}+y^{2})/4)$ (see \cite{OW} for details).
\end{example}

\section{Holomorphic $\infty$-harmonic maps}

In this section, we study $\infty$-harmonicity of holomorphic maps $
\mathbb{C}^{m}\longrightarrow \mathbb{C}^{n}$. Let
$(z_{1},\ldots,z_{m})\in \mathbb{C}^{m}$ and
$(w_{1},\ldots,w_{n})\in \mathbb{C}^{n}$  with
$z_{j}=x_{j}-iy_{j}\;j=1,\ldots, m$ and
$w_{\alpha}=u_{\alpha}-iv_{\alpha}\;\alpha=1,\ldots, n$. Then, a map
$\varphi: \mathbb{C}^{m}\longrightarrow \mathbb{C}^{n}$,
$\varphi(z_{1},\ldots,z_{m})=(\varphi^{1},\ldots, \varphi^{n})$ is
associated to a map $\varphi: \mathbb{R}^{2m}\longrightarrow
\mathbb{R}^{2n}$ with
$\varphi(x_{1},\ldots,x_{m},y_{1},\ldots,y_{m})=(u_{1},\ldots,u_{n},v_{1},\ldots,v_{n})$.
We write the map as $\varphi(X+iY)=\phi(X,Y)+i\psi(X,Y)$, where
$X=(x_{1},\ldots,x_{m}), Y=(y_{1},\ldots,y_{m})\in \mathbb{R}^{m}$
and the maps $\phi(X,Y)=(u^{1}(X,Y),\ldots, u^{n}(X,Y))$ and
$\psi(X,Y)=(v^{1}(X,Y),\ldots, v^{n}(X,Y))$ are called the real and
imaginary parts of $\varphi$. We have
\begin{theorem}\label{PE00}
A holomorphic  map $\varphi: \mathbb{C}^{m}\longrightarrow
\mathbb{C}^{n}$ with $\varphi(X,+iY)=\phi(X,Y)+i\psi(X,Y)$ is
$\infty$-harmonic if and only if its real and imaginary parts
$\phi(X,Y)$ and $\psi(X,Y)$ are $\infty$-harmonic.
\end{theorem}

\begin{proof}
It is well known that $\varphi: \mathbb{C}^{m}\longrightarrow
\mathbb{C}^{n}$ is holomorphic if and only if
\begin{equation}
\frac{\partial u^{\alpha}}{\partial x_{j}}=\frac{\partial
v^{\alpha}}{\partial y_{j}},\;\;\; \frac{\partial
u^{\alpha}}{\partial y_{j}}=-\frac{\partial v^{\alpha}}{\partial
x_{j}};\;\,\;\;\;j=1,2,\ldots,m,\;\alpha =1,2,\ldots,n.
\end{equation}
We can easily check that
\begin{eqnarray*}
|\nabla u^{\alpha}|^{2}&=&|\nabla v^{\alpha}|^{2},\\
\left|{\rm d}\varphi\right|^{2}&=&\delta
^{ij}{\varphi}^{\alpha}_{i}{\varphi}^{\beta}_{j}\delta_{\alpha
\beta} =\sum_{\alpha=1}^{n}|\nabla u^{\alpha}|^{2}+|\nabla
v^{\alpha}|^{2}=2\sum_{\alpha=1}^{n}|\nabla
u^{\alpha}|^{2}=2\sum_{\alpha=1}^{n}|\nabla v^{\alpha}|^{2},\\\notag
\nabla\left|{\rm
d}\varphi\right|^{2}&=&2\sum_{\alpha=1}^{n}\nabla|\nabla
u^{\alpha}|^{2}=2\sum_{\alpha=1}^{n}\nabla|\nabla v
^{\alpha}|^{2}\\&=& 2\nabla|\nabla \phi|^{2}= 2\nabla|\nabla
\psi|^{2}.\notag
\end{eqnarray*}
Substitute these into the $\infty$-harmonic map Equation
(\ref{00E2}) we obtain that $\varphi$ is $\infty$-harmonic if and
only if
\begin{equation}\label{N302910}
\begin{array}{lll}
2g(\nabla\,\phi^{\alpha},\nabla|\nabla \phi|^{2})=0, \\
2g(\nabla\,\psi^{\alpha},\nabla|\nabla \psi|^{2})=0,
\;\;\;\;\;\;\;\alpha =1,2,\ldots, n,
\end{array}
\end{equation}
which gives the Theorem.
\end{proof}

\begin{remark}
Explicitly, the $\infty$-harmonic map equation can be written as:
\begin{equation}\label{N0321}
\left\{\begin{array}{rl} 2\Delta_{\infty}u_{1}+\langle\nabla
u_{1},\nabla |\nabla
u_{2}|^{2}\rangle+\ldots+\langle\nabla u_{1},\nabla |\nabla u_{n}|^{2}\rangle=0\\
\langle\nabla u_{2},\nabla |\nabla
u_{1}|^{2}\rangle+2\Delta_{\infty}u_{2}+\ldots+\langle\nabla
u_{2},\nabla
|\nabla u_{n}|^{2}\rangle=0\\
\ldots \ldots \ldots \;\;\;\; \ldots \ldots \;\;\; \; \ldots \ldots \ldots\\
 \langle\nabla u_{n},\nabla |\nabla
u_{1}|^{2}\rangle+\langle\nabla u_{n},\nabla
|\nabla u_{2}|^{2}\rangle+\ldots++2\Delta_{\infty}u_{n}=0\\
2\Delta_{\infty}v_{1}+\langle\nabla v_{1},\nabla |\nabla
v_{2}|^{2}\rangle+\ldots+\langle\nabla v_{1},\nabla |\nabla v_{n}|^{2}\rangle=0\\
\langle\nabla v_{2},\nabla |\nabla
v_{1}|^{2}\rangle+2\Delta_{\infty}v_{2}+\ldots+\langle\nabla
v_{2},\nabla
|\nabla v_{n}|^{2}\rangle=0\\
\ldots \ldots \ldots \;\;\;\;\ldots \ldots \ldots \;\;\;\; \ldots \ldots \ldots\\
 \langle\nabla v_{n},\nabla |\nabla
v_{1}|^{2}\rangle+\langle\nabla v_{n},\nabla
|\nabla v_{2}|^{2}\rangle+\ldots+2\Delta_{\infty}v_{n}=0.\\
\end{array}\right.
\end{equation}
\end{remark}

\begin{theorem}\label{PE0010}
Let $\varphi: \mathbb{C}^{m}\longrightarrow \mathbb{C}$ be a
nonconstant holomorphic map. Then, $\varphi$ is an $\infty$-harmonic
map  if and only if $\varphi$ is  a composition of an orthogonal
projection $\mathbb{C}^{m}\longrightarrow \mathbb{C}$ followed by a
homothety $\mathbb{C}\longrightarrow \mathbb{C}$, i.e.,
$\varphi(z_{1},\ldots,z_{m})=\lambda z_{i} +z_{0}$, where $\lambda
\in\mathbb{R}, z_{0}\in\mathbb{C}$ are constants.
\end{theorem}
\begin{proof}
Notice that a holomorphic map  $\varphi:
\mathbb{C}^{m}\longrightarrow \mathbb{C}$ is automatically a
horizontally weakly conformal harmonic map (see e.g., \cite{BW1}).
It follows from the relationship among tension, $p$-tension, and
$\infty$-tension fields of a map
\begin{equation}\label{phm}
\tau_p(\varphi)= {\left|{\rm d }\varphi
\right|}^{p-2}\tau_{2}(\varphi)+(p-2){\left|{\rm d} \varphi
\right|}^{p-4}\tau_{\infty}(\varphi)
\end{equation}
 that if $\varphi$ is also an
$\infty$-harmonic map, then it must be $p$-harmonic for any $p$. In
this case $\varphi$ is a $p$-harmonic morphism (being horizontally
weakly conformal and $p$-harmonic map) for any $p$. By a theorem in
\cite{BL}, $\varphi$ must be a horizontally homothetic. Now the
classification of horizontally homothetic maps between Euclidean
spaces \cite{Ou1} applies to give the required results.
\end{proof}


\begin{thebibliography}{99}
\bibitem[Ar1]{Ar1} G. Aronsson, {\em Extension of functions satisfying Lipschitz conditions}, Ark. Mat. 6
(1967) 551--561.
\bibitem[Ar2]{Ar2} G. Aronsson, {\em On the partial differential equation
$u_{x}^{2}u_{xx} +2u_{x}u_{y}u_{xy}+u_{y}^{2}u_{yy}=0$}. Ark. Mat. 7
(1968) 395--425.
\bibitem[ACJ]{ACJ} G. Aronsson, M. Crandall, and P. Juutinen, {\em A tour of the theory of absolutely minimizing functions}
Bull. Amer. Math. Soc. (N.S.) 41 (2004), no. 4, 439--505.
\bibitem[BW]{BW1} P. Baird and  J. C. Wood, {\em Harmonic morphisms between Riemannian manifolds}, London Math. Soc. Monogr.
(N.S.) No. 29, Oxford Univ. Press (2003).
\bibitem[BB]{BB} G. Barles and J. Busca, {\em Existence and comparison results for fully nonlinear degenerate
elliptic equations without zeroth-order term}, Comm. Partial
Differential Equations 26 (2001), no. 11-12, 2323--2337.
\bibitem[Ba]{Ba} E. N. Barron, {\em Viscosity solutions and analysis in $L^\infty$}, in Nonlinear Analysis,
Differential Equations and Control (ed. by Clarke and Stern), Kluwer
Academic Publishers, 1999, 1-60.
\bibitem[BLW1]{BLW1} E. N. Barron, R. Jensen, and C. Y. Wang, {\em Lower Semicontinuity of $L^\infty$ functionals},
\bibitem[BLW2]{BLW2} E. N. Barron, R. Jensen, and C. Y. Wang, {\em Euler equations and absolute minimizers
of $L^\infty$ functionals}, Arch. Ration. Mech. Anal. 157 (2001),
no. 4, 255--283.
\bibitem[BEJ]{BEJ} E. N. Barron, L. C. Evans, and R. Jensen, The infinity Laplacian, Aronsson's equation and their generalizations, Preprint.
\bibitem[Bh]{Bh} T. Bhattacharya, {\em A note on non-negative singular infinity-harmonic functions in
he half-space}. Rev. Mat. Complut. 18 (2005), no. 2, 377--385.
\bibitem[BL]{BL} J. M. Burel and  E. Loubeau, {\em $p$-harmonic morphisms: the $1< p < 2$ case and a non-trivial example},
Contemp. Math. 308(2002), 21-37 .
\bibitem[CMS]{CMS} V. Caselles, J. -M. Morel, and C. Sbert, {\em An axiomatic approach to image interpolation}, IEEE Trans. Image Process.
7 (1998), no. 3, 376--386.
\bibitem[CEPB]{CEPB} G. Cong, M. Esser, B. Parvin, and G. Bebis, {\em Shape
metamorphism using $p$-Laplacian equation}, Proceedings of the 17th
International Conference on Pattern Recognition, (2004), Vol. 4,
15-18.
\bibitem[CE]{CE} M. G. Crandall and L. C. Evans, {\em A remark on infinity harmonic functions}. Proceedings of
the USA-Chile Workshop on Nonlinear Analysis (Vi$\tilde{{\rm n}}$a
del Mar-Valparaiso, 2000), 123--129.
\bibitem[CEG]{CEG} M. G. Crandall, L. C. Evans, and R. F. Gariepy, {\em Optimal Lipschitz extensions and the infinity Laplacian}.
Calc. Var. Partial Differential Equations 13 (2001), no. 2,
123--139.
\bibitem[CIL]{CIL} M. G. Crandall, H. Ishii, and P.-L. Lions, {\em User's guide to viscosity solutions of second order partial differential equations},
 Bull. Amer. Math. Soc. (N.S.) 27 (1992), no. 1, 1--67.
\bibitem[CY]{CY} M. G. Crandall and J. Zhang, {\em Another way to say harmonic}, Trans. Amer. Math. Soc. 355 (2003), 241-263.
\bibitem[EG]{EG} L. C. Evans and W. Gangbo, {\em Differential equations methods for the Monge-Kantorovich
mass transfer problem}, Mem. Amer. Math. Soc. 137 (1999), no. 653.
\bibitem[EY]{EY} L. C. Evans and Y. Yu, {\em Various properties of
solutions of the infinity -Laplace equation}, Preprint.
\bibitem[J]{J} R. Jensen, {\em Uniqueness of Lipschitz extensions: minimizing the sup norm of the gradient},
Arch. Rational Mech. Anal. 123 (1993), no. 1, 51--74.
\bibitem[JK]{JK} P. Juutinen and B. Kawohl, {\em On the evolution
governed by the infinity Laplacian}, Preprint.
\bibitem[JLM1]{JLM1} P. Juutinen, P. Lindqvist, and J. Manfredi,  {\em The infinity Laplacian: examples and observations}.
Papers on analysis, 207--217, Rep. Univ. Jyv$\ddot{{\rm
a}}$skyl$\ddot{{\rm a}}$ Dep. Math. Stat., 83, Univ. Jyv$\ddot{{\rm
a}}$skyl$\ddot{{\rm a}}$, Jyv$\ddot{{\rm a}}$skyl$\ddot{{\rm a}}$,
2001.
\bibitem[JLM2]{JLM2}P. Juutinen, P.
Lindqvist, and J. Manfredi, {\em The $\infty$ -eigenvalue problem},
Arch. Ration. Mech. Anal. 148 (1999), no. 2, 89--105.
\bibitem[LM1]{LM1} P. Lindqvist and J. Manfredi, {\em The Harnack inequality for
$\infty$ -harmonic functions}, Elec- tron. J. Differential Equations
(1995), No. 04, approx. 5 pp.
\bibitem[LM2]{LM2} P. Lindqvist and J. Manfredi, {\em Note on $\infty$ -superharmonic
functions}, Rev. Mat. Univ. Complut. Madrid 10 (1997), no. 2,
471--480.
\bibitem[Ob]{Ob} A. M. Oberman, {\em A convergent difference scheme for the infinity Laplacian:
construction of absolutely minimizing Lipschitz extensions}, Math.
Comp. 74 (2005), no. 251, 1217--1230 (electronic).
\bibitem[Ou1]{Ou1} Y. -L. Ou, {\em $p$-Harmonic morphisms, minimal
foliations, and rigidity of metrics}, J. Geom. Phys. 52 (2004), no.
4, 365--381.
\bibitem[Ou2]{Ou} Y. -L. Ou, {\em $p$-Harmonic functions and the minimal
graph equation in a Riemannain manifold}, Illinois J. Math. 49
(2005), no. 3, 911--927.
\bibitem[OTW]{OW} Y. -L. Ou, T. Troutman, and F. Wilhelm {\em $\infty$-harmonic maps and morphisms between
Riemannian manifolds}, preprint, 2007.
\bibitem[Sa]{Sa}  G. Sapiro, {\em Geometric partial differential equations and image analysis}, Cambridge University Press, Cambridge, 2001.
\bibitem[TP]{TP} M. Tenenbaum and H. Pollard, {\em Ordinary
differential equations}, Harper and Row, Publishers, New York, 1963.
\bibitem[Wa] [Wa] Ze-Ping Wang, Linear $\infty$-harmonic maps
between Riemannian manifolds, Preprint, 2007.
\bibitem[Zh]{Zh} R. Zhang, {\em $\infty$-Harmonic maps between semi-Euclidean
spaces}, Guangxi Sciences, (3) 2007, to appear.
\end{thebibliography}
\end{document}